\documentclass{amsart}
\usepackage{amsmath,amssymb,amsthm,amscd}
\newtheorem{formula}{}[section]
\newtheorem{proposition}[formula]{Proposition}
\newtheorem{corollary}[formula]{Corollary}
\newtheorem{lemma}[formula]{Lemma}
\newtheorem{theorem}[formula]{Theorem}
\theoremstyle{definition}
\newtheorem{definition}[formula]{Definition}
\newtheorem{example}[formula]{Example}
\theoremstyle{remark}
\newtheorem*{remark}{Remark}

\begin{document}

\title{Cohomology of graded Lie algebras of maximal class}
\author{Alice Fialowski}
\address{E\"otv\"os Lor\'and University, Department of Applied Analysis,
H-1117 Budapest, P\'azm\'any P. s\'et\'any, Hungary}
\email{fialowsk@cs.elte.hu}
\author{Dmitri Millionschikov}
\address{Moscow State University, Department of Mathematics and Mechanics,
Leninskie gory, 119899 Moscow, Russia}
\email{million@mech.math.msu.su}
\date{October 25, 2004}
\thanks{The research of the authors was partially supported by grants
OTKA T043641, T043034, RFBR 02-01-00659, "Russian Scientific Schools"
2185.2003.1. The final version
of the paper was completed during the stay of the second author at the
E\"otv\"os Lor\'and University in Budapest}
\begin{abstract}
It was shown by A. Fialowski that an arbitrary
infinite-dimensional ${\mathbb N}$-graded
"filiform type"  Lie algebra
$\mathfrak{g}{=}\bigoplus_{i{=}1}^{\infty}\mathfrak{g}_{i}$
with one-dimensional homogeneous components
$\mathfrak{g}_{i}$ such that
$[\mathfrak{g}_{1}, \mathfrak{g}_{i}]=\mathfrak{g}_{i+1},
\forall i \ge 2$ over a field of zero characteristic
is isomorphic to one (and only one)
Lie algebra from three given ones:
$\mathfrak{m}_0, \mathfrak{m}_2,  L_1$,
where the Lie algebras $\mathfrak{m}_0$ and $\mathfrak{m}_2$
are defined by their structure relations:
$\mathfrak{m}_0$: $[e_1,e_i]=e_{i+1}, \forall i \ge 2$ and
$\mathfrak{m}_2$: $[e_1, e_i ]=e_{i+1}, \forall i \ge 2$,
$[e_2, e_j ]=e_{j+2}, \forall j \ge 3$ and
$L_1$ is the "positive" part of the Witt algebra.

In the present article we compute the
cohomology
$H^*(\mathfrak{m}_0)$ and $H^*(\mathfrak{m}_2)$ with trivial coefficients,
give explicit formulas for their representative cocycles and describe
the multiplicative structure in the cohomology. Also we discuss
the relations with combinatorics and representation theory.
The cohomology $H^*(L_1)$ was calculated by L.~Goncharova in $1973$.
\end{abstract}
\date{}

\maketitle

\section*{Introduction}
${\mathbb N}$-graded Lie algebras are closely related to
nilpotent Lie algebras, for instance, a finite-dimensional
${\mathbb N}$-graded Lie algebra $\mathfrak{g}$ must to be nilpotent.
Infinite-dimensional ones are also called
residual nilpotent Lie algebras.
M.~Vergne studied in \cite{V} nilpotent
Lie algebras with the maximal possible nilindex
$s(\mathfrak{g})=\dim \mathfrak{g}{-}1$
(by nilindex $s(\mathfrak{g})$ we mean the
length of the descending central series
$\{C^i\mathfrak{g}\}$ of $\mathfrak{g}$). M.~Vergne called them
{\it filiform} Lie algebras. In her study the
${\mathbb N}$-graded filiform Lie algebra $\mathfrak{m}_0(n)$
played a special role.
This Lie algebra is defined by its basis $e_1,\dots,e_n$
and non-trivial commutator relations:
$[e_1,e_i]=e_{i+1}, i=2,\dots,n{-}1$.
Evidently, $\mathfrak{m}_0(n)$ is generated by $e_1$ and $e_2$.
Another example of ${\mathbb N}$-graded two-generated
filiform Lie algebra is $\mathfrak{m}_2(n)$:
$[e_1,e_i]=e_{i+1}, i=2,\dots,n{-}1, \:[e_2,e_i]=e_{i+2}, i=3,\dots,n{-}2$.

A.~Fialowski classified in \cite{Fial} all infinite-dimensional
${\mathbb N}$-graded two-generated Lie algebras
$\mathfrak{g}=\oplus_i\mathfrak{g}_i$
with one-dimensional homogeneous components $\mathfrak{g}_i$. In particular,
there are only three algebras in her list satisfying the "filiform
property": $[\mathfrak{g}_1,\mathfrak{g}_i]=\mathfrak{g}_{i+1}, \forall i$.
They are $\mathfrak{m}_0, \mathfrak{m}_2, L_1$, where
$\mathfrak{m}_0, \mathfrak{m}_2$ denote infinite-dimensional analogues
of $\mathfrak{m}_0(n), \mathfrak{m}_2(n)$, respectively and
$L_1$ is the "positive" part of the Witt or Virasoro algebra. The
classification of finite-dimensional ${\mathbb N}$-graded filiform
Lie algebras over a field of zero characteristic was done in \cite{Mill2}.

A.~Shalev and E.~Zelmanov defined in \cite{ShZ} the {\it coclass}
(which might be infinity) of a finitely generated and
residually nilpotent Lie algebra $\mathfrak{g}$, in analogy with the case
of (pro-)$p$-groups, as
$cc(\mathfrak{g})=
\sum_{i\ge 1}(\dim (C^i\mathfrak{g}/C^{i+1}\mathfrak{g})-1)$.
Obviously the coclass of a filiform algebra is equal to one
and the same is true
for the infinite-dimensional algebras $\mathfrak{m}_0, \mathfrak{m}_2, L_1$.
Algebras of coclass $1$ are also called algebras
of {\it maximal class}. They are also {\it narrow} or {\it thin} Lie algebras
(A.~Shalev, A.~Caranti, M.~Newman, et al.).
Part of Fialowski's classification in \cite{Fial} can be reformulated
in the following way: {\it Up to an isomorphism there are only
three ${\mathbb N}$-graded Lie algebras of maximal class with one-dimensional
homogeneous components:
$\mathfrak{m}_0, \mathfrak{m}_2, L_1$}. The last statement was
rediscovered in \cite{ShZ}.

Algebras of maximal class are in the center
of attention these days both in zero and positive characteristic.
There are many open questions related to them.
One natural question is their cohomology which is the subject of the present
paper.

L.~Goncharova calculated in $1973$ \cite{G}
the Betti numbers $b^q(L_1)=\dim H^q(L_1)$. Her result implies,
as a corollary,
the celebrated Euler identity in combinatorics:
$$\prod_{j=1}^{\infty}(1{-}t^j)
=\sum_{k=0}^{\infty}({-}1)^k(t^{\frac{3k^2{-}k}{2}}{+}
t^{\frac{3k^2{+}k}{2}})
$$
(see \cite{Fu}).
The Betti numbers $\dim H^q(\mathfrak{m}_0(n))$ (finite-dimensional case)
were calculated in \cite{AS} (see also \cite{B}), however, there are no
explicit formulas for basic cocycles and no
description of the multiplicative structure of $H^*(\mathfrak{m}_0(n))$
was obtained.

We give a complete description of the cohomology
algebras $H^*(\mathfrak{m}_0)$ and $H^*(\mathfrak{m}_2)$.
The method we use is based on Dixmier's exact sequence in
Lie algebra cohomology \cite{D}. In our considerations we use combinatorics:
partitions and generating functions.

The paper is organized as follows. In Sections 1--2 we review
all necessary definitions and facts concerning Lie algebra cohomology
and ${\mathbb N}$-graded Lie algebras, in particular we recall
Dixmier's exact sequence in the cohomology \cite{D}.
We start our computations in Section 3 with the algebra
$H^*(\mathfrak{m}_0)$ (Theorem \ref{main_H_m_0}).
In Section 4 we discuss the relations of our results with
representations theory. It turns out that
the basic cocycles representing $H^q(\mathfrak{m}_0)$
are at the same time the highest weight vectors (primitive vectors)
of the $q$-th exterior power $\Lambda^q(V(\lambda))$ of the irreducible
one-dimensional $\mathfrak{sl}(2,{\mathbb K})$-module $V(\lambda)$
for some $\lambda$ (Theorem \ref{sl_2-mod}).
In Section 5 we apply Dixmier's exact sequence
and the results of Section 3 to compute
$H^*(\mathfrak{m}_2)$ (Theorem \ref{main_H_m_2}).
Section 6 is devoted to finite-dimensional analogs of
the algebras considered above. Recall that all of them are filiform Lie
algebras. The Betti numbers $\dim H^q(\mathfrak{m}_0(n))$
are known \cite{AS}, \cite{B}. Some of $\dim H^q(L_1/L_{n+1})$
was found in \cite{Mill1}. The questions on the explicit formulas
for representing cocycles and the multiplicative structure are still
open for these algebras. At the end of the paper we consider
the characteristic $p$ analog of the algebras $\mathfrak{m}_0$
and $\mathfrak{m}_2$.
We briefly remark that other computational tools such as spectral
sequences
or the Hodge Laplacian of the differential $d$ lead to the same
algebraic and combinatorial problems that are solved in the present paper.

\section{Lie algebra cohomology and Dixmier's exact sequence}
Let us consider the cochain complex of a Lie algebra $\mathfrak{g}$
over a field ${\mathbb K}$ of zero characteristic:
$$
\begin{CD}
\mathbb K @>{d_0{=}0}>>
\mathfrak{g}^*=C^1(\mathfrak{g}) @>{d_1}>>
C^2 (\mathfrak{g}) @>{d_2}>>
\dots @>{d_{p-1}}>>C^{p} (\mathfrak{g}) @>{d_p}>> \dots
\end{CD}
$$
where $C^p(\mathfrak{g})$ denotes the vector space of continuous
skew-symmetric $p$-linear forms on $\mathfrak{g}$, and
the differential $d_p$ is defined by:
$$
d_pc(X_1, \dots, X_{p{+}1})=
\sum_{1{\le}i<j{\le}p{+}1}({-}1)^{i{+}j{-}1}
c([X_i,X_j], X_1,\dots, \hat X_j,\dots, \hat X_j, \dots, X_{p{+}1}).
$$
The differential
$d_1: \mathfrak{g}^* \rightarrow \Lambda^2 (\mathfrak{g}^*)$
coincides with the dual mapping of the Lie bracket
$[ \, , ]: \Lambda^2 \mathfrak{g} \to \mathfrak{g}$ and
$$
d(\rho \wedge \eta)=d\rho \wedge \eta+(-1)^{deg\rho} \rho \wedge d\eta,
\; \forall \rho, \eta \in \Lambda^{*} (\mathfrak{g}^*).
$$

\begin{definition}
\label{cohomology}
The cohomology of $(C^{*}(\mathfrak{g}), d)$ is
called the cohomology (with trivial coefficients) of the Lie algebra
$\mathfrak{g}$ and is denoted by $H^*(\mathfrak{g})$.
\end{definition}

Let $\mathfrak{b}$ be an
ideal of codimension $1$ in $\mathfrak{g}$.
One can choose an element $X\in \mathfrak{g}$ such that
$X \notin \mathfrak{b}$.
This element determines a $1$-form $\omega$ by the following condition:
$$\omega(X)=1,\;\: \omega(Y)=0\;\;{\rm for}\;\;  Y \in \mathfrak{b}.$$
As $\mathfrak{b}$ is an ideal in $\mathfrak{g}$, the $1$-form $\omega
\in {\mathfrak{g}}^*$ is closed and $\mathfrak{b}$ is an invariant
subspace of $adX$.

\begin{theorem}[Dixmier ~\cite{D}]
There exists a long exact sequence of Lie algebra cohomology:
\begin{equation}
\dots
\stackrel{(adX^*)_{i{-}1}}{\longrightarrow}
H^{i{-}1}(\mathfrak{b})
\stackrel{\omega \wedge}{\longrightarrow}
H^{i}(\mathfrak{g})
\stackrel{r_i}{\rightarrow} H^{i}(\mathfrak{b})
\stackrel{(adX^*)_i}{\longrightarrow} H^{i}(\mathfrak{b})
\to \dots \end{equation}
where

1) the homomorphism
$r_i: H^{i}(\mathfrak{g}) \to H^{i}(\mathfrak{b})$
is the restriction homomorphism;

2) $\omega \wedge: H^{*{-}1}(\mathfrak{b})
\to H^{*}(\mathfrak{g})$
is induced by the multiplication
$\omega \wedge: \Lambda^{*{-}1}(\mathfrak{b}^*)
\to \Lambda^*(\mathfrak{g}^*)$;

3) the homomorphisms
$(adX^*)_i: H^{i}(\mathfrak{b}) \to H^{i}(\mathfrak{b})$
are induced by the derivation $adX^*_i$ of degree
zero of $\Lambda^*(\mathfrak{b}^*)$ ($adX^*(a\wedge b)=
adX^*a\wedge b+ a\wedge adX^*b$ for all $a, b \in
\Lambda^*(\mathfrak{b}^*)$) that extends the dual mapping
$adX^* : \mathfrak{b}^* \to \mathfrak{b}^*$.
The $0$-derivation $adX^*$ commutes
with $d$ and we will denote the corresponding mapping in the cohomology
by the same symbol, more precisely, $(adX^*)_i$ denotes the mapping
$H^{i}(\mathfrak{b}) \to H^{i}(\mathfrak{b})$.
\end{theorem}
Dixmier's theorem \cite{D} is very important for our computations,
so it appears to be useful to recall his proof.
In fact, Dixmier considered the
cohomology of a nilpotent Lie algebra $\mathfrak{g}$ with coefficients
in arbitrary $\mathfrak{g}$-module $V$, but we restrict ourselves
to the case of scalar coefficients.
\begin{proof}
Each form $f \in \Lambda^*(\mathfrak{g})$
can be decomposed
as $f=\omega \wedge f' + f''$,
where $ f' \in \Lambda^{*{-}1}(\mathfrak{b}^*)$ and
$ f'' \in \Lambda^{*}(\mathfrak{b}^*)$.
One can write a short exact sequence of algebraic complexes
$$
0 \longrightarrow \Lambda^{*{-}1}(\mathfrak{b}^*)
\stackrel{\omega \wedge}{\longrightarrow} \Lambda^{*}(\mathfrak{g}^*)
\longrightarrow \Lambda^{*}(\mathfrak{b}^*)
\longrightarrow 0
$$
where $\Lambda^{*{-}1}(\mathfrak{b}^*)$
is taken with the differential $-d$ as
$d(\omega \wedge c)=d_{\lambda \omega}(\omega \wedge c)=
-\omega \wedge dc$.

This short exact sequence of algebraic complexes gives us
a long exact sequence in the cohomology.
To see this, we have to define the homomorphism
$H^{q}(\mathfrak{b}) \to H^{q}(\mathfrak{b})$
in the long exact sequence.

First of all let us introduce a new mapping
$$
\Lambda^{*}(\mathfrak{g}^*)
\to \Lambda^{*-1}(\mathfrak{b}^*),
f \in \Lambda^{*}(\mathfrak{g}^*)\to f_X \in
\Lambda^{*-1}(\mathfrak{b}^*),$$ where
$f_X(X_1,\dots,X_q)=f(X,X_1,\dots,X_q)$.
That means $f_X= f'$ if $f=\omega \wedge f' +f''$.

Then an obvious formula holds:
\begin{multline}
(d f)_X(X_1, \dots, X_{q{+}1})=
\sum_{1{\le}j{\le}q{+}1}({-}1)^{j}
f(adX(X_j), X_1,\dots, \hat X_j, \dots, X_{q{+}1})+\\
+\sum_{1{\le}i{<}j{\le}q{+}1}({-}1)^{i{+}j}
f([X_i,X_j],X,X_1, \dots, \hat X_i, \dots, \hat X_j, \dots, X_{q{+}1})=\\
=((adX^*)_q(f)+d(f_X))(X_1, \dots, X_{q{+}1}).
\end{multline}
Hence the homomorphism
$H^{q}(\mathfrak{b}) \to H^{q}(\mathfrak{b})$,
$[f] \to [(df)']$ of the long exact sequence coincides with the
homomorphism induced by $(adX^*)_q$.
\end{proof}
\section{${\mathbb N}$-graded Lie algebras}
\begin{definition}
A Lie algebra $\mathfrak{g}$ is called $\mathbb N$-graded,
if it is decomposed to the direct sum of subspaces such that
$$\mathfrak{g}=\oplus_{i} \mathfrak{g}_{i}, \; i \in
{\mathbb N},
\quad \quad [\mathfrak{g}_{i}, \mathfrak{g}_{j}]
\subset \mathfrak{g}_{i+j}, \;
\forall \: i, j \in {\mathbb N}.
$$
\end{definition}
\begin{example}
The Lie algebra $\mathfrak{m}_0$ is defined
by its infinite basis $e_1, e_2, \dots, e_n, \dots $
with commutator relations:
$$ \label{m_0}
[e_1,e_i]=e_{i+1}, \; \forall \; i \ge 2.$$
\end{example}
\begin{remark}
We always omit the trivial commutator
relations $[e_i,e_j]=0$ in the definitions of Lie algebras.
\end{remark}
\begin{example}
The Lie algebra $\mathfrak{m}_2$ is defined by its
infinite basis $e_1, e_2, \dots, e_n, \dots $
and commutator relations:
$$
[e_1, e_i ]=e_{i+1}, \quad \forall \; i \ge 2; \quad \quad
[e_2, e_j ]=e_{j+2}, \quad \forall \; j \ge 3.
$$
\end{example}
Consider now the algebra of polynomial vector fields
on the real line ${\mathbb R}^1$.

\begin{example}
Let us define the algebra $L_k$ as the Lie algebra
of polynomial vector fields
on the real line ${\mathbb R}^1$ with zero in
$x=0$ of order not less then $k+1$.
\end{example}
The algebra $L_k$ can be defined by its basis and the commutator relations are
given by the natural commutator of vector fields
$$e_i=x^{i+1}\frac{d}{dx}, \; i \in {\mathbb N},\; i \ge k, \quad \quad
[e_i,e_j]= (j-i)e_{i{+}j}, \; \forall \;i,j \in {\mathbb N}.$$
$L_k$ is a subalgebra of the Witt algebra $W$, where $W$ is
spanned by all $e_i=x^{i+1}\frac{d}{dx}, \; i \in {\mathbb Z}$.
Hence $L_1$ can be regarded as a "positive" part of the Witt
algebra.

The algebras $\mathfrak{m}_0, \mathfrak{m}_2, L_1$ are
${\mathbb N}$-graded Lie algebras generated by two elements $e_1, e_2$.
\begin{theorem}[follows from A.~Fialowski \cite{Fial}]
\label{fialowski}
Let $\mathfrak{g}=\bigoplus_{i=1}^{\infty}\mathfrak{g}_{i}$
be a ${\mathbb N}$-graded Lie algebra such that:
\begin{equation}
\dim \mathfrak{g}_{i}= 1, \; i \ge 1; \quad
[\mathfrak{g}_{1}, \mathfrak{g}_{i}]=\mathfrak{g}_{i+1},
\; \forall i \ge 2.
\end{equation}
Then $\mathfrak{g}$ is isomorphic to one (and only one)
Lie algebra from the three given ones:
$$\mathfrak{m}_0, \; \mathfrak{m}_2, \; L_1.$$
\end{theorem}
\begin{remark}
All  these Lie algebras have exactly two defining relations in
degree $5$ and $7$.
\end{remark}
The ideals $C^{k}\mathfrak{g}$ of the descending central series
of a Lie algebra $\mathfrak{g}$ determine
a decreasing filtration $C$ of $\mathfrak{g}$, i.e.
$[C^{k}\mathfrak{g},C^{l}\mathfrak{g}]
\subset C^{k+l}\mathfrak{g}, \: k,l \ge 1$.
One can consider the associated ${\mathbb N}$-graded
Lie algebra ${\rm gr}_C\mathfrak{g}$:
$${\rm gr}_C\mathfrak{g}=\oplus_{i\ge 1} ({\rm gr}_C\mathfrak{g})_i=
\oplus_{i\ge 1} C^i\mathfrak{g}/C^{i+1}\mathfrak{g}.$$
We have the following obvious isomorphisms of ${\mathbb N}$-graded
Lie algebras:
$$ {\rm gr}_C\mathfrak{m}_2 \cong
{\rm gr}_C L_1 \cong {\rm gr}_C\mathfrak{m}_0 \cong
\mathfrak{m}_0.$$
\begin{remark}
The Lie algebra $\mathfrak{m}_0$ has at least two
different ${\mathbb N}$-gradings: one of them we have already
discussed, and the second one (sometimes called natural) is defined by
$$\mathfrak{m}_0=\oplus_{i\ge 1}(\mathfrak{m}_0)_i, \;
(\mathfrak{m}_0)_1=\langle e_1, e_2 \rangle, \:
(\mathfrak{m}_0)_i=\langle e_{i+1} \rangle, \: i \ge 2.$$
The first grading is more convenient in cohomological
computations and we will use only this one.
\end{remark}

\begin{theorem}[Vergne \cite{V}]
\label{V_ne1}
Let $\mathfrak{g}=\oplus_{i\ge 1}\mathfrak{g}_{i}$ be a
${\mathbb N}$-graded Lie algebra such that
$$
\dim \mathfrak{g}_{1}=2, \: \dim \mathfrak{g}_{i}=1, \: i \ge 2; \;
[\mathfrak{g}_1,\mathfrak{g}_i]=\mathfrak{g}_{i+1}, \: i \ge 1.
$$
Then $\mathfrak{g}$ is isomorphic to $\mathfrak{m}_0$.
\end{theorem}

A Lie algebra $\mathfrak{g}$ in Definition \ref{cohomology}
is assumed to be topological
and the spaces $C^{*}(\mathfrak{g})$ can have
rather complicated nature, but here we will consider
only ${\mathbb N}$-graded Lie algebras
$\mathfrak{g}=\oplus_{i}\mathfrak{g}_{i}=
\oplus_{i}\langle e_{i} \rangle$.
Let us denote by $e^i, e^i(e_j)=\delta^i_j$ the corresponding dual $1$-forms.
In this situation we assume $C^{*}(\mathfrak{g})$ to be a vector space
of formal series of elements from $\Lambda^{*}(e^1, e^2, \dots)$
(we will also use the notation $\Lambda^{*}(\mathfrak{g}^*)$
for this space).

One can define a second grading in the
cochain complex $(C^*(\mathfrak{g}),d)=(\Lambda^{*}(\mathfrak{g}^*),d)$:
$\Lambda^{*}(\mathfrak{g}^*) =
\hat {\bigoplus_{k}} \Lambda^*_{k}(\mathfrak{g})$, where
a finite-dimensional subspace $\Lambda^q_{k}(\mathfrak{g})$ is spanned by
$q$-forms $e^{i_1} {\wedge} \dots {\wedge}
e^{i_q}, i_1<\dots <i_q$ such that
$i_1{+}\dots{+}i_q=k$. The symbol $\hat {\bigoplus_k}$
means the completed direct sum.

The second grading is compatible with the differential
$d$ and with the exterior product:
$$d \Lambda^q_{k} (\mathfrak{g})
\subset \Lambda^{q{+}1}_{k} (\mathfrak{g}), \qquad
\Lambda^{q}_{k}(\mathfrak{g}) \wedge
\Lambda^{p}_{l}(\mathfrak{g}) \subset
\Lambda^{q{+}p}_{k+l}(\mathfrak{g}).
$$
The exterior product in $\Lambda^*(\mathfrak{g})$
induces a structure of a
bigraded algebra in the cohomology
$H^*(\mathfrak{g})$:
$$
H^{q}_{k} (\mathfrak{g}) \wedge
H^{p}_{l} (\mathfrak{g}) \to
H^{q{+}p}_{k+l} (\mathfrak{g}).
$$
\begin{theorem}[Goncharova \cite{G}]
The Betti numbers $b^q(L_1)={\rm dim} H^q(L_1)=2$ for every
$q \ge 1$, more precisely
$$ ~~~b_k^q(L_1)={\rm dim} H_k^q(L_1)=
\left\{\begin{array}{r}
   1, \hspace{0.6em}{\rm if}~k=\frac{3q^2 \pm q}{2}, \\
   0,\hspace{1.76em}{\rm otherwise.}\\
   \end{array} \right . \hspace{3.3em} $$
\end{theorem}
The numbers
$\frac{3q^2 \pm q}{2}$ are so called Euler pentagonal numbers.
A sum of two arbitrary pentagonal numbers is not a pentagonal number,
hence the algebra
$H^*(L_1)$ has a trivial multiplication.
One can consider the Euler characteristic
$$\chi_k(L_1)=\sum_{q}({-}1)^q{\rm dim}C_k^q(L_1)=
\sum_{q}({-}1)^qb_k^q(L_1)$$
of each subcomplex $C^*_k(L_1), k \ge 0$, and
associate the corresponding generating function
$\chi_{L_1}(t)=\sum_{k=0}\chi_k(L_1)t^k$:
$$\prod_{j=1}^{\infty}(1{-}t^j){=}\sum_{k,q{\ge}0}({-}1)^q
{\rm dim}C_k^q(L_1)t^k{=}
\sum_{k,q{\ge}0}({-}1)^qb_k^q(L_1)t^k{=}
\sum_{k=0}^{\infty}({-}1)^k(t^{\frac{3k^2{-}k}{2}}{+}
t^{\frac{3k^2{+}k}{2}}).
$$
This equality of two expressions for the Euler characteristic
proves the celebrated Euler identity in combinatorics
(see \cite{Fu}, \cite{H} for details).
\begin{remark}
The ${\mathbb N}$-graded Lie algebras
$L_1$, $\mathfrak{m}_0$, $\mathfrak{m}_2$
have the same cochain spaces (but different differentials), hence
the cohomology $H^*(\mathfrak{m}_0)$ and $H^*(\mathfrak{m}_2)$
(that we are going to compute in the present article) also satisfy
the "Euler property":
$$\sum_{k=0}^{\infty}\sum_{q\ge 0}({-}1)^qb_k^q(\mathfrak{m}_0)t^k=
\sum_{k=0}^{\infty}\sum_{q\ge 0}({-}1)^qb_k^q(\mathfrak{m}_2)t^k=
\sum_{q=0}^{\infty}({-}1)^q(t^{\frac{3q^2{-}q}{2}}{+}
t^{\frac{3q^2{+}q}{2}}).
$$
\end{remark}

\section{The cohomology $H^*(\mathfrak{m}_0)$}
\label{H^*(m_0)}
Now we are going to apply Dixmier's exact sequence
for the situation
$$\mathfrak{g}=\mathfrak{m}_0, \; \mathfrak{b}= {\rm Span}
(e_2, e_3, \dots, e_k,\dots), \; X=e_1.$$

Obviously the ideal $\mathfrak{b}= {\rm Span}(e_2, e_3, \dots)$
is an abelian Lie algebra and its cohomology is the exterior algebra
$\Lambda^*(e_2, e_3, \dots)$.

Let us denote by $D_1$ the operator $ade_1^*: \Lambda^*(e_2, e_3, \dots)
\to \Lambda^*(e_2, e_3, \dots)$.
It can be defined explicitly by
\begin{equation}
\begin{split}
D_1(e^2)=0, \; D_1(e^i)= e^{i-1}, \; \forall i\ge 3,\\
D_1(\xi\wedge \eta)=D_1(\xi)\wedge \eta +\xi\wedge D_1(\eta),
\; \: \forall \xi, \eta \in \Lambda^*(e_2, e_3, \dots).
\end{split}
\end{equation}

\begin{lemma}
\label{surj}
The operator $D_1$ is surjective.
\end{lemma}
\begin{proof}
Let us define the operator $D_{-1}: \Lambda^*(e^2,e^3,\dots) \to
\Lambda^*(e^2,e^3,\dots)$,
\begin{equation}
\begin{split}
\label{D_{-1}}
D_{-1}e^i=e^{i+1}, \;
D_{-1}(\xi{\wedge} e^i)=
\sum_{l\ge 0}(-1)^l D_{1}^{l}(\xi){\wedge} e^{i+1+l},
\end{split}
\end{equation}
where $i\ge 2$ and $\xi$ is an arbitrary form in
$\Lambda^*(e^2,\dots,e^{i-1})$.
The sum in the definition (\ref{D_{-1}}) of $D_{-1}$ is always finite
because $D_1^l$ decreases the second grading by $l$.
For instance,
$D_{-1}(e^i\wedge e^k)=\sum_{l=0}^{i-2} ({-}1)^l e^{i-l}{\wedge} e^{k+l+1}$.

The operator $D_{-1}$ is right inverse to $D_1$, as
$D_1D_{-1}=Id$ on $\Lambda^*(e^2,e^3,\dots)$.

In fact
$D_1D_{-1}(e^i)=e^i, \; i \ge 2$ and for arbitrary
$\xi \in \Lambda^*(e^2,\dots,e^{i-1})$ we have
$$D_1D_{-1}(\xi\wedge e^i)=\sum_{l\ge 0}
(-1)^l D_{1}^{l+1}(\xi)\wedge e^{i+l+1}+ \sum_{l\ge 0}
(-1)^l D_{1}^{l}(\xi)\wedge e^{i+l}=\xi\wedge e^{i}.$$
\end{proof}
One can write the formula for
$D_{-1}(e^{i_1} {\wedge} \dots {\wedge} e^{i_q}{\wedge} e^{i_q})=D_{-1}(0)$:
\begin{equation}
\label{m_0cocycles}
\omega(e^{i_1}{\wedge}\dots {\wedge} e^{i_q} {\wedge}
e^{i_q{+}1})=
\sum\limits_{l\ge 0}(-1)^l D_1^l(e^{i_1}\wedge \dots
\wedge e^{i_q})\wedge e^{i_q+1+l}.
\end{equation}
This sum is also always finite and determines
a homogeneous closed $(q{+}1)$-form
of the second grading $i_1{+}{\dots}{+}i_{q{-}1}{+}2i_{q}{+}1$.

Let us consider the restriction
$(D_1)_k^{q+1}$ of $D_1$ on $(q+1)$-forms of the second grading $k$.
\begin{lemma} \label{KerD_1} Let $q\ge 1$.
Then ${\rm Ker}(D_1)_k^{q+1}$ is spanned by
$$
\omega(e^{i_1}{\wedge}\dots {\wedge} e^{i_q} {\wedge}
e^{i_q{+}1}), \:2\le i_1{<}{\dots}{<}i_q,
\:i_1{+}{\dots}{+}i_{q-1}{+}2i_q{+}1=k.
$$
\end{lemma}
\begin{proof}
Obviously, if these forms exist, they are linearly independent.
How many they are? Let us consider
the inclusion $A^{q+1}_{k-1}:\Lambda_{k-1}^{q+1}(e^2,e^3,\dots)
\to \Lambda_{k}^{q+1}(e^2,e^3,\dots)$ defined on basic monomials
by the shift of the last superscript:
$$
A^{q+1}_{k-1}(\xi{\wedge}e^i{\wedge}e^j)=
\xi{\wedge}e^i{\wedge}e^{j{+}1}.
$$
It follows immediately that ${\rm Im}A^{q+1}_{k-1}$ is spanned by monomials
$e^{j_1}{\wedge}\dots{\wedge}e^{j_{q}}{\wedge}e^{j_{q+1}}$,
such that $j_{q+1}-j_q\ge 2$. Hence $\dim {\rm Im} A^{q+1}_{k-1}=
\dim\Lambda_{k-1}^{q+1}(e^2,e^3,\dots)$ and
$$
{\rm Im}A^{q+1}_{k-1}\oplus
{\rm Span}
(\{\omega(e^{i_1}{\wedge}\dots{\wedge}e^{i_q}{\wedge}e^{i_q+1}),
i_1{+}{\dots}{+}2i_q{+}1=k \})=
\Lambda_{k}^{q+1}(e^2,e^3,\dots).
$$
On the another hand, $(D_1)_k^{q+1}$ is surjective
and hence
$$\dim {\rm Ker}(D_1)_k^{q+1}=\dim \Lambda_{k}^{q+1}(e^2,e^3,\dots)-
\dim \Lambda_{k-1}^{q+1}(e^2,e^3,\dots)$$
which completes the proof.
\end{proof}

We will think of $\omega$ as a linear map defined on
the subspace of $\Lambda^{q+1}(e^2,e^3,\dots)$, spanned by the monomials
$e^{i_1}{\wedge}\dots{\wedge}e^{i_q}{\wedge}e^{i_q+1}$:
$$
\omega(\sum_{i_1\dots i_{q}} \alpha_{i_1\dots i_{q}}
e^{i_1}{\wedge}\dots{\wedge}e^{i_q}{\wedge}e^{i_q+1}):=
\sum_{i_1\dots i_{q+1}} \alpha_{i_1\dots i_{q+1}}
\omega(e^{i_1}{\wedge}\dots{\wedge}e^{i_q}{\wedge}e^{i_q+1}).
$$
Let us denote by $P_q(k)$ the number of (unordered) partitions
of a positive integer $k$ into $q$ parts, i.e.
(see \cite{H}) $P_q(k)$ is the number
of solutions in positive integers $x_i$ of
$$
k=x_1+\dots+x_q, \; 1\le x_1 \le x_2 \le \dots \le x_q.
$$
Denote by $V_q(k)$ the number of partitions of a positive integer
$k$ into $q$ {\it distinct} parts, i.e. $V_q(k)$ is the number
of solutions in positive integers $y_i$ of
$$
k=y_1+\dots+y_q, \; 1< y_1 < y_2 < \dots < y_q.
$$
\begin{lemma}
\label{partitions}
${\rm dim} \Lambda_k^q(e^2, e^3, \dots)=
V_q\left(k{-}q\right)=
P_q\left(k{-}\frac{q(q{+}1)}{2}\right)$.
\end{lemma}
\begin{proof}
It consists of two standard tricks in combinatorics.
Let us consider

1) a one-to-one correspondence between the set
of partitions of a positive integer $N$ exactly into
$q$ parts and the standard basis of
$\Lambda_{N{+}\frac{q(q{+}1)}{2}}^q(e^2, e^3, \dots)$:
\begin{multline*}
(x_1, x_2, x_3, \dots, x_q) \to e^{x_1{+}1}{\wedge}
e^{x_2{+}2}{\wedge}e^{x_3{+}3}{\wedge}\dots {\wedge}e^{x_q{+}q},\\
N=x_1+\dots+x_q, \quad 1\le x_1 \le x_2 \le \dots \le x_q;
\end{multline*}

2) an analogous bijection between the partitions of $N$ into
$q$ distinct parts and the standard basis of
$\Lambda_{N{+}q}^q(e^2, e^3, \dots)$:
\begin{multline*}
(y_1, y_2, y_3, \dots, y_q) \to e^{y_1{+}1}{\wedge}
e^{y_2{+}1}{\wedge}e^{y_3{+}1}{\wedge}\dots {\wedge}e^{y_q{+}1},\\
N=y_1+\dots+y_q, \quad 1\le y_1 < y_2 < \dots < y_q.
\end{multline*}
\end{proof}
\begin{theorem}
\label{main_H_m_0}
The bigraded cohomology algebra
$H^*(\mathfrak{m}_0)=\oplus_{k,q} H^q_k(\mathfrak{m}_0)$
is spanned by the cohomology classes of the following homogeneous cocycles:
\begin{equation}
\label{sumHm_0}
\begin{split}
e^1, \: e^2, \:
\omega(e^{i_1}{\wedge}\dots
{\wedge} e^{i_q}{\wedge} e^{i_q{+}1})=
\sum\limits_{l\ge 0}({-}1)^l (ad e_1^*)^l
(e^{i_1}{\wedge} e^{i_2}{\wedge} \dots
{\wedge} e^{i_q}){\wedge} e^{i_q{+}1{+}l}, \\
\end{split}
\end{equation}
where
$q \ge 1, \; 2\le i_1 {<}i_2{<}{\dots} {<}i_q$,
in particular
$${\rm dim} H_{k{+}\frac{q(q{+}1)}{2}}^q(\mathfrak{m}_0)=
P_q(k)-P_q(k{-}1).$$
The multiplicative structure is defined by
\begin{multline}
\label{Hm_0multiplication}
[e^1] {\wedge} \omega(\xi{\wedge}e^i{\wedge}e^{i+1})=0,\;
[e^2] {\wedge} \omega(\xi{\wedge}e^i{\wedge}e^{i+1})
=\omega(e^2{\wedge}\xi{\wedge}e^i{\wedge}e^{i+1}),\\
\omega(\xi{\wedge}e^i{\wedge}e^{i+1})
{\wedge} \omega(\eta{\wedge}e^j{\wedge}e^{j+1})
=\sum_{l=0}^{j-i+1}(-1)^l\omega((ade_1^*)^l(\xi{\wedge}e^i){\wedge}e^{i+1+l}
{\wedge}\eta{\wedge}e^j{\wedge}e^{j+1})+\\
+(-1)^{j{-}i{+}deg\eta}\sum_{s \ge 1}
\omega((ade_1^*)^{j{-}i{-}1{+}s}(\xi{\wedge}e^i){\wedge}
(ade_1^*)^s(\eta{\wedge}e^j){\wedge}e^{j{+}s}
{\wedge}e^{j{+}s{+}1})+\\
+(-1)^{j{-}i{+}deg\eta {+}1}\sum_{s \ge 1}
\omega((ade_1^*)^{j{-}i{+}1{+}s}(\xi{\wedge}e^i){\wedge}
(ade_1^*)^s(\eta{\wedge}e^j){\wedge}e^{j{+}s{+}1}
{\wedge}e^{j{+}s{+}2}),
\end{multline}
where $i \le j$, $\xi$ and $\eta$ are arbitrary homogeneous forms in
$\Lambda^*(e^2,\dots,e^{i{-}1})$ and $\Lambda^*(e^2,\dots,e^{j{-}1})$,
respectively.
\end{theorem}
\begin{proof}
It follows from Lemma \ref{surj} that
in our case Dixmier's sequence is equivalent to the following exact
sequences:
\begin{equation*}
\begin{split}
0\stackrel{}{\longrightarrow}
H^{0}(\mathfrak{b})={\mathbb K}
\stackrel{e^1 {\wedge}}{\longrightarrow}
H^1(\mathfrak{m}_0)
\stackrel{r_1}{\longrightarrow} {\langle} e^2 \rangle
\to 0;\\
0\stackrel{}{\longrightarrow}
H^q(\mathfrak{m}_0)
\stackrel{r_q}{\longrightarrow} {\rm Ker}(ade_1^*)_q
\to 0, \; q \ge 2.
\end{split}
\end{equation*}
Here ${\rm Ker}(ade_1^*)_q={\rm Ker}(D_1)^q$ and the dimensions
$\dim {\rm Ker}(D_1)^q_k$ were found in Lemma \ref{KerD_1}
and the final formula for
$\dim H^q_{k{+}\frac{q(q{+}1)}{2}}(\mathfrak{m}_0)$
follows from Lemma \ref{partitions}.

Let us remark that the sum (\ref{sumHm_0}) is always finite. Namely, the
maximal value of the superscript $l$ is equal
to $l_{max}=i_1{+}\dots{+}i_q-\frac{q(q+3)}{2}$ and the corresponding
summand in
(\ref{sumHm_0}) is
$(-1)^{l_{max}}\alpha(i_1,\dots, i_q)e^{2}{\wedge}e^{3}{\wedge}\dots
{\wedge} e^{q{+}2}{\wedge} e^{i_q{+}1{+}l_{max}}$
(see Example \ref{3cocycle} below).

We will obtain the formulas for the multiplicative structure of
 $H^*(\mathfrak{m}_0)$ by using the explicit expressions (\ref{m_0cocycles})
for basic cocycles.
An arbitrary homogeneous $q$-cocycle
$\omega(e^{i_1}{\wedge}\dots {\wedge}
e^{i_{q{-}2}}{\wedge}e^r{\wedge}e^{r{+}1})$ is completely determined
by its leading term
$e^{i_1}{\wedge}\dots {\wedge}e^{i_{q{-}2}}{\wedge}e^r{\wedge}e^{r{+}1}$
 -- the unique monomial
$\alpha e^{i_1}{\wedge}{\dots}e^{i_{q{-}1}}{\wedge}e^{i_q},
\:i_1<\dots < i_{q{-}1}<i_q$ in its decomposition
such that $i_q-i_{q{-}1}=1$.
Hence we only have to consider the summands with this property in the
right part of the formula
$$
\omega(\xi{\wedge}e^i{\wedge}e^{i{+}1}){\wedge}
\omega(\eta{\wedge}e^j{\wedge}e^{j{+}1}){=}
\sum_{l,k\ge 0}({-}1)^{k{+}l}D_1^l(\xi{\wedge}e^i){\wedge}
e^{i{+}1{+}l}{\wedge}D_1^k(\eta{\wedge}e^j){\wedge}e^{j{+}1{+}k}.
$$
They are of the following three kinds:

{\parindent 0pt
\medskip
1) \ $({-}1)^{l}D_1^l(\xi{\wedge}e^i){\wedge}
e^{i{+}l{+}1}{\wedge}\eta{\wedge}e^j{\wedge}e^{j{+}1}\quad
{\rm for} \quad i+1+l \le j+2;$
\smallskip

2) \  $({-}1)^{j{-}i{-}1}D_1^{j{-}i{-}1{+}k}
(\xi{\wedge}e^i){\wedge}
e^{j{+}k}{\wedge}D_1^k(\eta{\wedge}e^j){\wedge}e^{j{+}1{+}k} \:
\quad{\rm for}\quad k\ge 1 ;$
\smallskip

3) \ $({-}1)^{j{-}i{+}1}D_1^{j{-}i{+}1{+}k}
(\xi{\wedge}e^i){\wedge}
e^{j{+}2{+}k}{\wedge}D_1^k(\eta{\wedge}e^j){\wedge}e^{j{+}1{+}k} \:
\quad{\rm for}\quad k\ge 1.$

\medskip

Taking $\omega$ of their sum we get formula (\ref{Hm_0multiplication})
for the multiplication in $H^*(\mathfrak{m}_0)$ (recall that $D_1$ denotes
the operator $ade_1^*$).
}
\end{proof}

\begin{example}
We choose the following basis of $H^2(\mathfrak{m}_0)$:
$$
e^2{\wedge} e^3, e^3{\wedge} e^4 -e^2{\wedge} e^5, \dots,
\omega(e^j{\wedge}e^{j+1})=
\sum_{l=0}^{j-2} ({-}1)^l e^{j-l} {\wedge} e^{j+1+l}, \dots
$$
>From this it is clear that
$$  ~~~{\rm dim} H_k^2(\mathfrak{m}_0)= \left\{\begin{array}{r}
   1, \; k=2j+1 \ge 5, \\
   0, \; {\rm otherwise} \quad \quad \; \\
   \end{array} \right.$$
and for $2 \le i < j$ we have
\begin{multline*}
\omega(e^2{\wedge}e^{3}){\wedge}\omega(e^j{\wedge}e^{j{+}1})=
\omega(e^2{\wedge}e^{3}{\wedge}e^j{\wedge}e^{j{+}1}), \\
\omega(e^i{\wedge}e^{i+1})
{\wedge} \omega(e^j{\wedge}e^{j+1})
{=}{\sum_{l{=}0}^{j{-}i{+}1}{(}{-}1{)}^l}\omega{(}e^{i{-}l}{\wedge}e^{i+1+l}
{\wedge}e^j{\wedge}e^{j+1}{)}+\\
+({-}1)^{j{-}i}\sum_{s=1}^{2i{-}j{-}1}\omega(e^{2i{-}j{+}1{-}s}{\wedge}
e^{j{-}s}{\wedge}e^{j{+}s}{\wedge}e^{j{+}s{+}1})+\\
+({-}1)^{j{-}i{+}1}
\sum_{s=1}^{2i{-}j{-}3}
\omega (e^{2i{-}j{-}1{-}s}{\wedge}
e^{j{-}s}{\wedge}e^{j{+}s{+}1}
{\wedge}e^{j{+}s{+}2}).
\end{multline*}
\end{example}
\begin{corollary} $H_k^3(\mathfrak{m}_0)=0$ if $k < 9$,

and for $k\ge 9$ we
have
$$  ~~~{\rm dim} H_k^3(\mathfrak{m}_0)= \left\{\begin{array}{r}
   l-1, \; k=6l{+}r, \: r=0{,}1{,}2{,}4,\\
   l, \quad \quad \quad k=6l{+}r, \: r=3{,}5.\\
   \end{array} \right . \hspace{3.3em} $$
\end{corollary}
This follows from the formulas for $P_3(k)$ given in \cite{H}.
Using other remarks on $P_q(k)$ in \cite{H},
one can show that ${\rm dim} H_k^q(\mathfrak{m}_0)$
is a polynomial of degree $q-2$ in $k$ with  leading term
$\frac{k^{q{-}2}}{(q{-}2)!q!}$ and the other coefficients depend
on residue $k$ modulo $q!$.
\begin{example}
\label{3cocycle}
\begin{multline*}
\omega(e^5{\wedge} e^6 {\wedge} e^7)=
e^5{\wedge} e^6 {\wedge}e^7-e^4{\wedge} e^6 {\wedge} e^8+(e^3{\wedge} e^6+
e^4 {\wedge} e^5){\wedge} e^9-\\
-(e^2{\wedge}e^6+2e^3 {\wedge}e^5){\wedge}e^{10}
+(3e^2{\wedge} e^5+2e^3{\wedge} e^4){\wedge} e^{11}-
5e^2{\wedge} e^4{\wedge} e^{12}+5e^2{\wedge} e^3{\wedge} e^{13}.
\end{multline*}
\end{example}

Now we can finish this section with the formula for the generating
function for Betti numbers $b^q_k(\mathfrak{m}_0)$:
$$
\sum_{k=0}^{\infty}\sum_{q=0}^k b^q_k(\mathfrak{m}_0)t^kx^q=
t(1+x)+(1-t)\prod_{j=2}^{\infty}(1+xt^j),
$$
and verify the Euler property:
\begin{multline*}
\sum_{k=0}^{\infty}\sum_{q\ge 0}(-1)b^q_k(\mathfrak{m}_0)t^k=
1-t-t^2+\sum_{k\ge 3}t^k\sum_{q\ge 2}({-}1)^q(V_q(k{-}q){-}V_q(k{-}q{-}1))=\\
=1{-}t{-}t^2{+}(\prod_{j=2}^{\infty}(1{-}t^j)-1+t^2)-
t(\prod_{j=2}^{\infty}(1{-}t^j)-1)=(1{-}t)\prod_{j=2}^{\infty}(1{-}t^j)=
\prod_{j=1}^{\infty}(1{-}t^j).
\end{multline*}

\section{The cohomology $H(\mathfrak{m}_0)$ and highest weight representations of
$\mathfrak{sl}(2,{\mathbb K})$}
This section was influenced by Appendix B in
Bordemann's article \cite{B}.

Let $X, Y, H$ denote the standard basis of
$\mathfrak{sl}(2,{\mathbb K})$:
$$[X,Y]=H, \; [H,X]=2X, \; [H,Y]=-2Y.$$
One can define an infinite-dimensional
$\mathfrak{sl}(2,{\mathbb K})$-module $V(\lambda)$ by its basis
$\{f_i, \: i \ge 0\}$ and the well-known classical formulas (\cite{S}):
\begin{equation}
\begin{split}
Hf_i=&(\lambda-2i)f_i,\\
Yf_i=&(i+1)f_{i+1},\\
Xf_i=&(\lambda-i+1)f_{i-1},
\end{split}
\end{equation}
where we set $f_{-1}=0$ and $\lambda \in {\mathbb K}$.
The $\mathfrak{sl}(2,{\mathbb K})$-module $V(\lambda)$
is generated by its highest weight vector (primitive vector \cite{S})
$f_0$:
$Hf_0=\lambda f_0$,
$Xf_0=0$, $f_i=Y^if_0/i!$. The module $V(\lambda)$ is
irreducible if and only if $\lambda \notin {\mathbb N}, \lambda \ne 0$.
Sometimes the module $V(\lambda)$ is called the standard cyclic module or
the Verma module.

Consider now the $q$-th exterior power $\Lambda^q(V(\lambda))$.
It is decomposed to the direct sum of its weight subspaces
$(\Lambda^q(V(\lambda)))^{\alpha}=\{v \in \Lambda^q(V(\lambda)),
H(v)=\alpha v\}$:
$$
\Lambda^q(V(\lambda))=\oplus_k (\Lambda^q(V(\lambda)))^{\lambda q-2k}=
\oplus_k \Lambda^q_k(V(\lambda)),
$$
where  $\Lambda^q_k(V(\lambda))$ is spanned by monomials
$f_{i_1}{\wedge}\dots{\wedge}f_{i_q}$ such that
$i_1+\dots+i_q=k$. Obviously
$H(f_{i_1}{\wedge}\dots{\wedge}f_{i_q})=(\lambda q -2\sum_{l=1}^qi_l)
f_{i_1}{\wedge}\dots{\wedge}f_{i_q}$.

Now we take $\lambda\notin {\mathbb N}, \lambda \ne 0$ and consider
a new basis of $V(\lambda)$:
$$
\label{tilde_f}
{\tilde f}_0=f_0,\: {\tilde f}_i=f_i/\prod_{l=1}^i(\lambda-l+1),\:
i\ge 1, \; X{\tilde f}_i={\tilde f}_{i-1}, \: i \ge 0.
$$
In fact, we can reformulate Theorem \ref{main_H_m_0} as follows.
\begin{theorem} \label{sl_2-mod}
Let $q \ge 2$ and $V(\lambda)$ be an infinite-dimensional
irreducible $\mathfrak{sl}(2,{\mathbb K})$-module, then
its $q$-th exterior power $\Lambda^q(V(\lambda))$ is reducible.
The subspace of its primitive vectors of weight $\lambda q - 2k$ has
dimension $P_q(k-\frac{(q-3)q}{2})-P_q(k-\frac{(q-3)q}{2}-1)$ and
we can take the following basis:
\begin{multline}
\omega({\tilde f}_{i_1}{\wedge}\dots
{\wedge}{\tilde f}_{i_{q-2}}{\wedge}{\tilde f}_{i}{\wedge}{\tilde f}_{i{+}1})
=\sum\limits_{l\ge 0}({-}1)^l X^l
({\tilde f}_{i_1}{\wedge}\dots{\wedge}{\tilde f}_{i_{q-2}}{\wedge}
{\tilde f}_i){\wedge} {\tilde f}_{i{+}1{+}l}, \\
0 \le i_1 < i_2 < \dots < i_{q-2}< i, \; \sum_{l=1}^{q-2}i_l+2i+1=k.
\end{multline}
That means the module
$\Lambda^q(V(\lambda))$ is decomposed to the direct sum of its highest weight
submodules $\Lambda^q(V(\lambda))=\oplus V_{i_1\dots i_{q-2},i}$,
where the highest weight vector of $V_{i_1\dots i_{q-2},i}$ is
$\omega({\tilde f}_{i_1}{\wedge}\dots{\wedge}{\tilde f}_{i_{q-2}}
{\wedge}{\tilde f}_{i}{\wedge}{\tilde f}_{i{+}1})$.
\end{theorem}
\begin{corollary}
Let $\lambda \notin {\mathbb Q}$ and $q \ge 2$, then each
submodule $V_{i_1\dots i_{q-2},i}$
in the decomposition of $\Lambda^q(V)$ considered above
is irreducible (and hence isomorphic to $V(\lambda)$).
\end{corollary}
\begin{proof}
The operator $D_1: {\rm Span}(e^2,e^3,\dots) \to
{\rm Span}(e^2,e^3,\dots)$ considered above
plays the role of the operator $X$
in an infinite-dimensional irreducible
$\mathfrak{sl}(2,{\mathbb K})$-module
$V(\lambda)={\rm Span}(e^2,e^3,\dots)$. We just rescale
${\tilde f}_i=e^{i+2},i \ge 0$, where $\{{\tilde f}_i, \: i \ge 0\}$
is the basis (\ref{tilde_f}) of $V(\lambda)$.
The subspace of primitive elements of weight $\lambda q - k$
coincides by definition with $({\rm Ker}D_1)^q_{k+2q}$ and
the theorem follows from Lemma~\ref{KerD_1}.
\end{proof}

\section{The cohomology $H^*(\mathfrak{m}_2)$}

Let us set in Dixmier's exact sequence
$$\mathfrak{g}=\mathfrak{m}_2, \; \mathfrak{b}= {\rm Span}
(e_1,e_3,\dots), \; X=e_2.$$
The ideal $\mathfrak{b}$ is isomorphic to
$\mathfrak{m}_0$, by shifting $e_1=e_1',\: e_i=e_{i-1}', i\ge3$.
Hence it follows from the results of the previous section
that the cohomology $H^*(\mathfrak{b})$
is spanned by the cohomology classes of
$$e^1, \:e^3, \:
\omega_{\mathfrak{b}}(e^{i_1}{\wedge}\dots {\wedge} e^{i_p} {\wedge}
e^{i_p{+}1}), \; p\ge 1,\: 3\le i_1 <\dots < i_p,$$
where in the formula
$\omega_{\mathfrak{b}}(\xi{\wedge}e^{i} {\wedge}
e^{i+1})=\sum_{l\ge 0}D_1^l(\xi{\wedge}e^{i}){\wedge}e^{i+1+l},
\xi \in \Lambda^*(e^2,\dots,e^{i-1})$ and
we also assume that $D_1(e^3)=0$, for instance
$\omega_{\mathfrak{b}}(e^{3}{\wedge}e^4)=e^{3}{\wedge}e^4$.

The operator $D_2=ade_2^*: \Lambda^*(\mathfrak{b})
\to \Lambda^*(\mathfrak{b})$
inducing ${\mathcal D}_2: H^*(\mathfrak{b}) \to H^*(\mathfrak{b})$
can be defined by
\begin{equation}
\begin{split}
D_2(e^1)=D_2(e^4)=0,\; D_2(e^3)=e^1,\; D_2(e^i)= e^{i-2},
\: i\ge 5,\\
D_2(\xi\wedge \eta)=D_2(\xi)\wedge \eta +\xi\wedge D_2(\eta), \;
\forall \xi, \eta \in \Lambda^*(\mathfrak{b}).
\end{split}
\end{equation}

It is immediate that
\begin{equation}
\label{2-surj}
\begin{split}
{\mathcal D}_2(e^3)=e^1&, \quad {\mathcal D}_2(e^1)=0,\\
{\mathcal D}_2(\omega_{\mathfrak{b}}(e^3{\wedge} e^4))=0,& \;
{\mathcal D}_2(\omega_{\mathfrak{b}}(e^k{\wedge} e^{k+1}))=
-2\omega_{\mathfrak{b}}(e^{k-1}{\wedge} e^k), \;
k \ge 4.
\end{split}
\end{equation}
\begin{proposition} Let $3\le i_1<\dots<i_{p-1}<i$ and
$\xi=e^{i_1}{\wedge}\dots{\wedge}e^{i_{p-1}}$, then
$$
\label{mathcal D_2}
{\mathcal D}_2(\omega_{\mathfrak{b}}(\xi{\wedge}e^{i}{\wedge}e^{i+1}))=
\omega_{\mathfrak{b}}((D_2+D_1^2)(\xi){\wedge}e^{i}{\wedge}e^{i+1})
-2\omega_{\mathfrak{b}}(\xi{\wedge}e^{i-1}{\wedge}e^{i}).$$
\end{proposition}
\begin{proof}
First of all let us remark that in general $D_2$ and $D_1^2$ coincide
only on $\Lambda^1(\mathfrak{b})$:
$$D_1^2(e^i)=D_2(e^i)=e^{i-2}, \:i\ge 5, \;
D_1^2(\xi{\wedge}\eta)=D_1^2(\xi){\wedge}\eta+
2D_1(\xi){\wedge}D_1(\eta)+\xi{\wedge}D_1^2(\eta).$$
An arbitrary cocycle is completely determined by its terms
$\alpha e^{j_1}{\wedge}\dots{\wedge}e^{j_{p}}{\wedge}
e^{j_p+1}$. On the another hand, the operator $D_2$ decreases the
difference between the two last superscripts of some
monomial $e^{j_1}{\wedge}\dots{\wedge}e^{j_{p}}$ by two.
Denoting all "non-interesting" terms by dots, we obtain
an expression for
$D_2(\omega_{\mathfrak{b}}(\xi{\wedge}e^{i}{\wedge}e^{i+1}))$:
\begin{multline*}
D_2(\xi{\wedge}e^{i}{\wedge}e^{i+1}-
D_1(\xi{\wedge}e^{i}){\wedge}e^{i+2}+
D_1^2(\xi{\wedge}e^{i}){\wedge}e^{i+3}+\dots)=\\
=D_2(\xi){\wedge}e^{i}{\wedge}e^{i+1}+
\xi{\wedge}e^{i}{\wedge}D_2(e^{i+1})-\\
-\xi{\wedge}D_1(e^{i}){\wedge}D_2(e^{i+2})+
D_1^2(\xi){\wedge}e^{i}{\wedge}D_2(e^{i+3})+\dots=\\
=D_2(\xi){\wedge}e^{i}{\wedge}e^{i+1}+
\xi{\wedge}e^{i}{\wedge}e^{i-1}
-\xi{\wedge}e^{i-1}{\wedge}e^{i}+
D_1^2(\xi){\wedge}e^{i}{\wedge}e^{i+1}+\dots.
\end{multline*}
\end{proof}
\begin{lemma}
$
H^*(\mathfrak{b})={\rm Im}{\mathcal D}_2 \oplus \langle e^3 \rangle.
$
\end{lemma}
\begin{proof}
We have already seen (\ref{2-surj}) that ${\mathcal D}_2$ is surjective
on $H^2(\mathfrak{b})$. Now
for an arbitrary $\xi=e^{i_1}{\wedge}\dots{\wedge}e^{i_{q-1}}$
with $q \ge 2$ and $3\le i_1<\dots<i_{q-1}<i$, define an operator
${\mathcal D}_{-2}:H^2(\mathfrak{b}) \to H^2(\mathfrak{b})$ by the formula:
$$
\label{DD_{-2}}
{\mathcal D}_{-2}(\omega_{\mathfrak{b}}(\xi {\wedge}e^{i}{\wedge} e^{i+1}))=
\sum_{l\ge 0}\frac{1}{2^{l+1}}\omega_{\mathfrak{b}}\left((D_2+D_1^2)^l
(\xi) {\wedge}
e^{i+l}{\wedge} e^{i+1+l}\right).
$$
Like in the definition of $D_{-1}$, remark that the sum (\ref{DD_{-2}})
is always finite: the operator $(D_2+D_1^2)^l$ decreases the second
grading of an arbitrary homogeneous $\xi$ by $2l$.
Now using formula (\ref{mathcal D_2}) we obtain
${\mathcal D}_{2}{\mathcal D}_{-2}=-Id$.
\end{proof}
Let us consider the  restriction
$({\mathcal D}_2)_k^{q+1}={\mathcal D}_2:H^{q+1}_k(\mathfrak{b}) \to
H^{q+1}_{k-2}(\mathfrak{b})$.
\begin{corollary} \label{dimH_m_2}
Let $q\ge 1$. Then
$\dim {\rm Ker}({\mathcal D}_2)_k^{q+1}=
\dim H_{k}^{q+1}(\mathfrak{b})-
\dim H_{k-2}^{q+1}(\mathfrak{b})$.
\end{corollary}

\begin{lemma} \label{KerD_2} The space
${\rm Ker}{\mathcal D}_2$ is spanned by
$$e^1, \; \omega_{\mathfrak{b}}(e^3{\wedge}e^4)=e^3{\wedge}e^4, \;
\sum_{l\ge 1}\frac{1}{2^{l}}\omega_{\mathfrak{b}}\left((D_2+D_1^2)^l
(e^{i_1}{\wedge} {\dots} {\wedge} e^{i_q}) {\wedge}
e^{i_q{+}1{+}l}{\wedge} e^{i_q{+}2{+}l}\right),
$$
where $1\le q, \; 3\le i_1 <i_2<\dots <i_q$.
\end{lemma}
\begin{proof}
A mimic of the proof of Lemma \ref{KerD_1}.
\end{proof}

\begin{theorem}
\label{main_H_m_2}
The bigraded cohomology algebra
$H^*(\mathfrak{m}_2)=\oplus_{q,k} H^q_k(\mathfrak{m}_2)$ is spanned
by cohomology classes of the following homogeneous cocycles:
\begin{multline}
e^1, \; e^2, \;  e^2 \wedge e^3, \; e^3\wedge e^4-e^2\wedge e^5, \\
w_{i_1{,} {\dots}{,} i_q{,} i_q{+}1{,} i_q{+}2}=
\sum_{l{\ge}1}\frac{1}{2^l}\omega\left({(}ade_2^*{+}{(}ade_1^*{)}^2{)}^l
{(}e^{i_1}{\wedge}{\dots} {\wedge} e^{i_q}{)} {\wedge}
e^{i_q{+}1{+}l}{\wedge} e^{i_q{+}2{+}l}\right),
\end{multline}
where $1\le q, \; 3\le i_1 <i_2<\dots <i_q$, in particular
for $q \ge 3$,
$$
{\rm dim} H_{k{+}\frac{q(q{+}1)}{2}}^q(\mathfrak{m}_2)=
P_q(k)-P_q(k{-}1)-P_q(k{-}2)+
P_q(k{-}3).
$$
\end{theorem}
\begin{proof}
Dixmier's sequence is equivalent to the following exact
sequences:
\begin{multline*}
0\stackrel{}{\longrightarrow}
H^{0}(\mathfrak{b})={\mathbb K}
\stackrel{e^2 {\wedge}}{\longrightarrow}
H^1(\mathfrak{m}_0)
\stackrel{r_1}{\longrightarrow} {\langle} e^1 \rangle
\longrightarrow 0;\\
0\stackrel{}{\longrightarrow}\langle e^3 \rangle
\stackrel{e^2 {\wedge}}{\longrightarrow}
H^2(\mathfrak{m}_2)
\stackrel{r_2}{\longrightarrow}
\langle \omega_{\mathfrak{b}}(e^3{\wedge}e^4) \rangle
\longrightarrow 0,\\
0\stackrel{}{\longrightarrow}
H^q(\mathfrak{m}_2)
\stackrel{r_q}{\longrightarrow} {\rm Ker}(ade_2^*)_q
\longrightarrow 0, \; q \ge 3.
\end{multline*}
The equality ${\rm Ker}(ade_2^*)_q={\rm Ker}({\mathcal D}_2)^q$
was found in Lemma \ref{KerD_2}.
formulas for $\dim H^q_k(\mathfrak{m}_2)$
follows from Corollary \ref{dimH_m_2} and Theorem \ref{main_H_m_0}.
The last remark is that $\omega(\xi)$ represents the inverse image of
$\omega_{\mathfrak{b}}(\xi) \in {\rm Ker}(ade_2^*)_q$ with respect
to $r_q:H^q(\mathfrak{m}_2) \to H^q(\mathfrak{b})$, for instance
$r_q^{-1}(\omega_{\mathfrak{b}}(e^3{\wedge}e^4))=
r_q^{-1}(e^3{\wedge}e^4)=\omega(e^3{\wedge}e^4)=
e^3{\wedge}e^4-e^2{\wedge}e^5.$
\end{proof}
\begin{example}
\begin{multline*}
w_{5,6,7}=\omega(e^5{\wedge} e^6 {\wedge} e^7) +
\omega(e^3{\wedge} e^7{\wedge}e^8)
=e^5{\wedge} e^6 {\wedge}e^7+(e^3{\wedge} e^7{-}e^4{\wedge} e^6){\wedge} e^8+\\
{+}(e^4 {\wedge} e^5{-}e^2{\wedge} e^7){\wedge} e^9+
(e^2{\wedge}e^6{-}e^3 {\wedge}e^5){\wedge}e^{10}
+e^3{\wedge} e^4{\wedge} e^{11}-
e^2{\wedge} e^4{\wedge} e^{12}+e^2{\wedge} e^3{\wedge} e^{13}.
\end{multline*}
\end{example}
\begin{corollary}

1) The space $H^2(\mathfrak{m}_2)$ is two-dimensional
and it is spanned by the cohomology classes represented by
cocycles $e^2\wedge e^3$ and $e^3\wedge e^4-e^2\wedge e^5$
with second grading $5$ and $7$ respectively;\\
2) The space $H^3(\mathfrak{m}_2)$ is
infinite-dimensional and is spanned by
$$w_{k,k{+}1,k{+}2}=
\sum_{l\ge 0}\omega\left(e^{k-2l}{\wedge}
e^{k+1+l}{\wedge} e^{k+2+l}\right), \; k\ge 3.$$
Hence
$  ~~{\rm dim} H_q^3(\mathfrak{m}_2)= \left\{\begin{array}{r}
   1, \; q=3k{+}3 \ge 12, \\
   0, \; {\rm otherwise.} \quad \quad \; \\
   \end{array} \right . \hspace{3.3em} $
\end{corollary}

\begin{remark}
Again (see \cite{H}), one can show that $b^q_k(\mathfrak{m}_2)$
for $q \ge 3$ is a polynomial of degree $q-3$ with leading term
$\frac{2k^{q{-}3}}{(q{-}3)!q!}$ and other coefficients depending
on residue $k$ mod $q!$. We have the following
identity for the corresponding generating function:
$$
\sum_{k=0}^{\infty}\sum_{q=0}^k b^q_k(\mathfrak{m}_2)t^kx^q=
(1+x)(t+t^2-t^3+xt^5)+(1-t-t^2+t^3)\prod_{j=3}^{\infty}(1+xt^j).
$$
It is easy to verify that the Euler property of $b^q_k(\mathfrak{m}_2)$
is equivalent to the following obvious equality:
$(1{-}t{-}t^2{+}t^3)\prod_{j=3}^{\infty}(1{-}t^j)=
\prod_{j=1}^{\infty}(1{-}t^j)$.
\end{remark}

\section{Finite-dimensional analogs and their cohomology}
Let $\mathfrak{g}$ be a ${\mathbb N}$-graded Lie algebra,
then $\oplus_{i=n{+}1} \mathfrak{g}_i$ is an ideal and we
can consider the corresponding quotient Lie algebra
$\mathfrak{g}/\oplus_{i=n{+}1} \mathfrak{g}_i$ and denote it
by $\mathfrak{g}(n)$.
The quotient Lie algebras $\mathfrak{m}_0(n)$, $\mathfrak{m}_2(n)$ and
${\mathcal V}_n=L_1/L_{n+1}$ are $n$-dimensional nilpotent Lie algebras
with the same length $s(\mathfrak{g})=n-1$
of the descending central series $\{C^i\mathfrak{g}\}$.
In fact, $n-1$ is the maximum of $s(\mathfrak{g})$ in the
set of $n$-dimensional nilpotent Lie algebras.
Studying this class of nilpotent Lie algebras was initiated
by Vergne in \cite{V}.

Recall that
a nilpotent $n$-dimensional Lie algebra $\mathfrak{g}$ is  called
filiform Lie algebra, if its descending central series $\{C^i\mathfrak{g}\}$
has the length (nil-index) $s(\mathfrak{g})=n-1$.

The filiform Lie algebra $\mathfrak{m}_0(n)$ plays a special role in the
theory of filiform Lie algebras: an arbitrary $n$-dimensional
filiform Lie algebra can be obtained as a nilpotent deformation
of $\mathfrak{m}_0(n)$ (\cite{V}).
Vergne's explicit formulas for
basic cocycles $\Psi_{k,r}$ of the second cohomology
$H^2(\mathfrak{m}_0, \mathfrak{m}_0)$
with coefficients in the adjoint representation
are one of the main tools in the study of filiform Lie algebras
via the deformation theory. We can mention that
the cohomology groups $H^2(\mathfrak{m}_0)$ and
$H^2(\mathfrak{m}_2)$ were also found by M.~Vergne (\cite{V}).

Later all the Betti numbers of $\mathfrak{m}_0(n)$ were calculated
by M.~Bordemann in \cite{B}(Appendix B). We would like to discuss
his very elegant approach, in fact similar results are well known
in combinatorics (R.~Stanley's  and his school, see \cite{AS} for
references). By arguments equivalent to Dixmier's exact sequence,
he reduced the computation of $\dim H^q(\mathfrak{m}_0(n))$ to the
problem of finding ${\rm Ker}(D_1)_q$:
$$
\dim H^q(\mathfrak{m}_0(n))=\dim{\rm Ker}(D_1(n))_q+
\dim{\rm Ker}(D_1(n))_{q-1}.
$$
Let us take $\lambda=n-2$ and consider the irreducible
$(n-1)$-dimensional $\mathfrak{sl}(2,{\mathbb K})$-module $V(n-2)=
{\rm Span}(e^2,\dots,e^n)$, where $X=D_1$.
The dimension of ${\rm Ker}(D_1(n))_q$ is equal to
the number of irreducible $\mathfrak{sl}(2,{\mathbb K})$-modules
in the decomposition of $\Lambda^q(V(n-2))$ and the last number
is equal in its
turn to the dimension of the zero-eigenspace plus
the dimension of one-eigenspace of
$H:\Lambda^q(V(n{-}2)) {\to} \Lambda^q(V(n{-}2))$.
Now we rescale as before ${\tilde f}_i=e^{i+2}, \: i=0,\dots,n-2$.
$$H({\tilde f}_{i_1}{\wedge}\dots{\wedge}{\tilde f}_{i_q})=
(q\lambda-2\sum_{l=1}^qi_l)
{\tilde f}_{i_1}{\wedge}\dots{\wedge}{\tilde f}_{i_q}, \; \lambda=n-2.
$$
Hence $\dim{\rm Ker}(D_1(n))_q$ is equal to the number of solutions
of the equation
$$
\sum_{l=1}^qi_l=\left[\frac{q(n-2)}{2}\right], \;
0 \le i_1 <\dots <i_q\le n-2.
$$
In other words,
$\dim{\rm Ker}(D_1(n))_q=V_{q,n-1}(\left[\frac{qn}{2}\right])$,
where we denoted by $V_{q,n-1}(N)$ the number of partitions
of a positive integer $N$ into $q$ distinct summands
$i_1,\dots,i_q$ such that $1\le i_1 <\dots <i_q\le n{-}1$ and
$[x]$ stands for the integer part of $x \in {\mathbb Q}$.
We conclude that
\begin{equation}
\dim H^q(\mathfrak{m}_0(n))=V_{q,n-1}\left(\left[qn/2\right]\right)+
V_{q-1,n-1}\left(\left[(q-1)n/2\right]\right).
\end{equation}
It is possible to write some of the first formulas in terms of more
convenient
combinatorial functions:
\begin{equation*}
\begin{split}
\dim H^2(\mathfrak{m}_0(n))=\left[ \frac{n+1}{2}\right],\;
\dim H^3(\mathfrak{m}_0(n))=\left[ \binom{\frac{n+1}{2}}{2}+
\frac{1}{8}\right],\\
\dim H^4(\mathfrak{m}_0(n))=\left[ \frac{4}{3}\binom{\frac{n+1}{2}}{3}+
\frac{4n+13}{36}\right].
\end{split}
\end{equation*}
\begin{remark}

The last formulas give no information about the bigraded structure
of $H^*(\mathfrak{m}_0(n))$. Also we have no explicit formulas for
basic cocycles. Later Bordemann's results were generalized in
\cite{AS} to the case of an arbitrary  finite-dimensional
nilpotent Lie algebra with an abelian ideal of codimension one
(the preprint version of \cite{B} appeared earlier than
\cite{AS}).
\end{remark}

It was shown in \cite{Mill1} that
the case $\mathfrak{g}(n)$ can be more complicated then
its infinite-dimensional analog $\mathfrak{g}$.
The Betti numbers of $L_1/L_{n+1}$ stabilize as
$n \to \infty$  and $\dim H^q(L_1/L_{n+1})$
are equal to the Fibonacci numbers $F_{q+2}$ for sufficiently large $n>>q$.
For instance,
$$\dim H^2(L_1/L_{n+1})=3,\:\dim H^3(L_1/L_{n+1})=5,\:
\dim H^4(L_1/L_{n+1})=8,\:\dots.$$
However, the full description of $H^*(L_1/L_{n+1})$ is still
an open question. The same situation is with the cohomology
$H^*(\mathfrak{m}_2)$.

The classification of finite-dimensional ${\mathbb N}$-graded filiform Lie
algebras $\mathfrak{g}=\oplus\langle e_i \rangle$
with one-dimensional homogeneous components
(there is the one-parametric family $\mathfrak{g}_{\alpha}$
of non-isomorphic algebras
in each dimension $7\le \dim \mathfrak{g} \le 11$, see \cite{Mill2})
also shows that one can expect the difficulties in the
finite-dimensional case.

\section{Final remarks}
1) We can consider the Lie algebras $\mathfrak{m}_0$ and
$\mathfrak{m}_2$ not only over a field of zero characteristic but
over an arbitrary field of positive characteristic.
>From the proof of Theorem \ref{main_H_m_0} it follows that the statement
is also valid over an arbitrary field and Theorem \ref{main_H_m_2}
is valid over any field of non-even characteristic.
This remark appears important
because of applications of Lie algebras of maximal class in the
theory of (pro-)$p$-groups (see \cite{ShZ}.

2) Let $\mathfrak{g}=\oplus_{i\ge 1}\langle e_i \rangle$
be an ${\mathbb N}$-graded Lie algebra with one-dimensional
homogeneous components.
We can equip the (finite-dimensional)
space $\Lambda^q_k(\mathfrak{g}^*)$
with an euclidean scalar product $(,)$, such that the basic
monomials $e^{i_1}{\wedge}\dots{\wedge}e^{i_q}$ form
an orthonormal basis of $\Lambda^q_k(\mathfrak{g}^*)$:
$$
(e^{i_1}{\wedge}\dots{\wedge}e^{i_q},
e^{j_1}{\wedge}\dots{\wedge}e^{j_q})=
\delta^{i_1}_{j_1}\dots\delta^{i_q}_{j_q}.
$$
One of the main methods to compute the cohomology $H^*(\mathfrak{g})$
of the cochain complex $(\Lambda^*(\mathfrak{g}),d)$ is
to find the zero-eigenspace of the Hodge Laplacian $dd^*+d^*d$
(see \cite{Fu}).

\begin{proposition}
The Hodge Laplacian $dd^*+d^*d$ of the Lie algebra $\mathfrak{m}_0$
satisfies the following properties:
\begin{equation}
\begin{split}
(dd^*+d^*d)(e^1)=0,& \: \
(dd^*+d^*d)(e^1{\wedge}\xi)=e^1{\wedge}D_1D_1^*(\xi),\\
(dd^*+d^*d)(\eta)&=D_1^*D_1(\eta),\;
\quad\eta \notin e^1{\wedge}\Lambda^*(\mathfrak{g}^*)
\end{split}
\end{equation}
where $D_1^*=ade_1$ is the $0$-derivation of the exterior algebra
$\Lambda^*(\mathfrak{g})$ that extends the operator $ade_1:\mathfrak{m}_0 \to
\mathfrak{m}_0$.
\end{proposition}
Let us study the kernel of the Hodge Laplacian.
We have proved in Lemma \ref{surj} that
$D_1: \Lambda^*(e^2,e^3,\dots)\to \Lambda^*(e^2,e^3,\dots)$
is surjective, hence ${\rm Ker}D_1^*=0$, moreover it is easy to see that
$e^1 \notin {\rm Im} D_1$ and on the another hand,
${\rm Ker}D_1D_1^*=0$. Hence
${\rm Ker}(dd^*+d^*d) \subset e^1\oplus \Lambda^*(e^2,e^3,\dots)$
and the problem of determining the zero-eigenspace
of $dd^*+d^*d$ reduces to the computation
of ${\rm Ker}D_1^*D_1={\rm Ker}D_1$. The basic cocycles
$\omega(e^{i_1}{\wedge}\dots{\wedge} e^{i_q}{\wedge} e^{i_q{+}1})$
from Theorem \ref{main_H_m_0} together with
$e^1, e^2$ form the basis of the harmonic forms.

3) One can compute the cohomology of $\mathfrak{m}_0$ directly
without Dixmier's exact sequence, but for the Lie algebra
$\mathfrak{m}_2$ such a computation is not so easy. We chose
Dixmier's method (see also \cite{AS}) to demonstrate a general
principle which is applicable for other finite or infinite
dimensional Lie algebras as well. Dixmier remarked in \cite{D}
that his exact sequence in the cohomology is equivalent to some
cohomological spectral sequence.

Let $\mathfrak{g}=\oplus_i \mathfrak{g}_i$
be an arbitrary ${\mathbb N}$-graded Lie algebra
of the filiform type, i.e. $\dim \mathfrak{g}_i =1$ for all $i$ and
$[\mathfrak{g}_1,\mathfrak{g}_i]=\mathfrak{g}_{i+1}, \: i \ge 2$.
We showed earlier
that it is filtered by the ideals $\{C^i\mathfrak{g}\}$
of the descending central series and the corresponding
associated graded Lie algebra ${\rm gr}_C \mathfrak{g}$ is isomorphic
to $\mathfrak{m}_0$. It is easy to see that this filtration induces
a filtration in the cochain complex
$\Lambda^*(\mathfrak{g})$ and the first term $E_1$
of the corresponding spectral sequence $(E_r, d_r)$ is isomorphic to the
cohomology $H^*({\rm gr}_C \mathfrak{g})=H^*(\mathfrak{m}_0)$.
Moreover, the spectral sequence $(E_r, d_r)$ degenerates at the first term
$E_1$ in the case $\mathfrak{g}=\mathfrak{m}_0$,
and at the second term $E_2$ in the case $\mathfrak{g}=\mathfrak{m}_2$,
respectively. In both cases the differential $d_0$
is equal to $e^1{\wedge}D_1$, and $d_1=0$ for $\mathfrak{g}=\mathfrak{m}_0$,
while $d_1=e^2{\wedge}D_2$ for $\mathfrak{g}=\mathfrak{m}_2$.
The case of $\mathfrak{g}=L_1$ is much more complicated. We have infinite
number of non-trivial differentials $d_r$ of the spectral sequence $E_r$
and computing all of them would give a new proof of
Goncharova's theorem. This approach might be helpful to obtain
explicit formulas for Goncharova's cocycles
which have not yet been found so far.

\end{document}